\input amstex
\magnification=1200
\documentstyle{amsppt}
\nologo
\hoffset=-0.5pc \vsize=57.2truepc \hsize=37truepc
\spaceskip=.5em plus.25em minus.20em
\define\verti{\Cal V}
\define\frakv{\nu}
\define\lambd{\frak q}
\define\QQ{M}
\define\qq{b}
\define\guisteft{1}
\define\marleone{2}
\define\ortrathr{3}
\define\perrosou{4}
\define\schmah{5}
\define\sjamlerm{6}

\topmatter
\title Slices for lifted tangent and cotangent actions
\endtitle
\author Johannes Huebschmann
\endauthor
\affil Universit\'e des Sciences et Technologies de Lille
\\
U. F. R. de Math\'ematiques
\\
CNRS-UMR 8524
\\
F-59 655 VILLENEUVE D'ASCQ, France
\\
Johannes.Huebschmann\@math.univ-lille1.fr
\endaffil
\date{August 24, 2005}
\enddate
\abstract{Given a Lie group $G$, a $G$-manifold $\QQ$, and a point
$\qq$ of $\QQ$ with compact stabilizer, we construct slices for
the lifted tangent and cotangent actions at a pre-image of $\qq$
in terms of a slice for the $G$-action  on $\QQ$ at the point
$\qq$. We interpret the slice for the lifted cotangent action in
terms of a symplectic slice and in terms of a
Witt-Artin-decomposition.}
\endabstract

\address{\smallskip
\noindent USTL, UFR de Math\'ematiques, CNRS-UMR 8524
\newline\noindent
F-59 655 Villeneuve d'Ascq C\'edex, France
\newline\noindent
Johannes.Huebschmann\@math.univ-lille1.fr}
\endaddress
\endkeywords
\endtopmatter
\document
\beginsection Introduction

Let $G$ be a Lie group,  $\QQ$ a smooth manifold, suppose that
$\QQ$ is endowed with a smooth $G$-action, let $\qq$ be a point of
$\QQ$, and suppose that the stabilizer $G_{\qq}$ of $\qq$ is
compact. In this paper we will explore $G$-slices at a pre-image
of the point $\qq$ for the lifted $G$-actions on the tangent and
cotangent bundles of $\QQ$ in terms of a $G$-slice at $\qq$.

To this end, we endow $\QQ$ with a $G_{\qq}$-invariant Riemannian
metric. Then the tangent space $\roman T_{\qq} \QQ$ to $\QQ$ at
the point $\qq$ is an orthogonal $G_{\qq}$-representation, and the
orthogonal complement $S_{\qq}=(\frak g \qq)^{\bot}$ in $\roman
T_{\qq} \QQ$ of the tangent space $\frak g \qq$ to the $G$-orbit
$G\qq$ in $\QQ$ is well known to be an infinitesimal slice at
$\qq$ for the $G$-action on $\QQ$, that is, for a suitable
$G_{\qq}$-invariant ball $B_{\qq}$ in $S_{\qq}$ containing the
origin, the smooth map from $G \times _{G_{\qq}}B_{\qq}$ to $M$
which is given by the assignment to $(x,y)\in G \times B_{\qq}$ of
$x\cdot\roman{exp}_{\qq}(y)\in M$ is a $G$-equivariant
diffeomorphism onto a neighborhood of the $G$-orbit $G\qq$ of
$\qq$ in $\QQ$ in such a way that the zero section is identified
with $G\qq$. The  orthogonal representation of $G_{\qq}$ on
$S_{\qq}$ is then referred to as the {\it slice representation\/}
at the point $\qq$.  The total spaces $\roman T\QQ$ and $\roman
T^*\QQ$ of the tangent and cotangent bundles, respectively, of
$\QQ$ inherit smooth $G$-actions. In this paper we shall construct
infinitesimal slices for these lifted $G$-actions on $\roman T\QQ$
and $\roman T^*\QQ$ and in particular describe the resulting slice
representations entirely in terms of data that involve only the
base manifold $\QQ$ and the geometry of the group $G$. In
particular, the slice representation for the cotangent bundle case
given in Theorem 3.9 below implies the following.

\proclaim{Theorem} Let $\frak g \qq \ (\subseteq \roman
T_{\qq}(\QQ))$ be the tangent space to the $G$-orbit $G\qq$ in
$\QQ$ at the point $\qq$ of $\QQ$, let $S_{\qq}=(\frak g
\qq)^{\bot}\subseteq \roman T_{\qq}(\QQ)$ be the corresponding
infinitesimal slice at $\qq$ for the $G$-action on $\QQ$, and let
$\roman T^*_{\qq}(\QQ)= (\frak g \qq)^* \oplus S^*_{\qq}$ be the
resulting orthogonal decomposition of $\roman T^*_{\qq}(\QQ)$.
Moreover, let $\frak g_{\qq}p$ be the tangent space to the
$G_{\qq}$-orbit $G_{\qq}p$ (in $\roman T^*_{\qq}(\QQ)$) at the
point $p$ of $\roman T^*_{\qq}(\QQ)$, necessarily a linear
subspace of $S^*_{\qq}\cong\roman T_{\qq} S^*_{\qq} \subseteq
\roman T^*_{\qq}(\QQ)$, let $S^*_{\qq} = \frak g_{\qq}p \oplus
\sigma^*_{\qq} $ be the resulting orthogonal decomposition of
$S^*_{\qq}$, and let $(\frak g_{\qq}p)^*\subseteq \roman T_{\qq}
S_{\qq}\cong S_{\qq}$ be the vector space which is dual to $\frak
g_{\qq}p$ so that $S_{\qq} = (\frak g_{\qq}p)^* \oplus
\sigma_{\qq} $ is the corresponding orthogonal decomposition of
$S_{\qq}$. Then, at the point $p$ of $\roman T^*_{\qq}(\QQ)$, as
an orthogonal $G_p$-representation, a suitable infinitesimal slice
for the lifted $G$-action on the total space $\roman T^*\QQ$ of
the cotangent bundle of $\QQ$ is isomorphic to the orthogonal
$G_p$-representation
$$
(\frak g \qq)^*\oplus \sigma^*_{\qq} \oplus (\frak g_{\qq}p)^*
\oplus \sigma_{\qq}.
$$
\endproclaim

This description of the infinitesimal $G$-slice at the point $p$
of the total space $\roman T^*\QQ$ of the cotangent bundle of
$\QQ$ plainly involves  only data phrased in terms of the base
$\QQ$. Notice that, in particular, $p$ is a point of $\roman
T^*_{\qq}(\QQ)$, that $\frak g_{\qq}p$ is the tangent space to the
$G_{\qq}$-orbit $G_{\qq}p$ (in $\roman T^*_{\qq}(\QQ)$) at the
point $p$ of $\roman T^*_{\qq}(\QQ)$, and that $(\frak
g_{\qq}p)^*\subseteq \roman T_{\qq} S_{\qq}\cong S_{\qq}$ is the
vector space which is dual to $\frak g_{\qq}p$. From this
description of the infinitesimal $G$-slice at $p$, it is
straightforward to derive a {\it symplectic slice\/} and a
corresponding {\it Witt-Artin decomposition\/} at $p$; for details
and comments related with the significance of the Witt-Artin
decomposition, see \cite\ortrathr.

At the present stage, we can already explain the symplectic slice
and the Witt-Artin decomposition at the point $p$ without too much
trouble, leaving some of the details to Section 4 below; to this
end, let $\mu \colon \roman T^*\QQ@>>> \frak g^*$ be the cotangent
bundle momentum mapping for the $G$-action on $\roman T^*\QQ$, let
$\frakv = \mu(p)\in \frak g^*$, and let $\Cal O_{\frakv} = G\frakv
\subseteq \frak g^*$, the coadjoint orbit generated by $\frakv$.
Let $G_{\frakv}\subseteq G$ be the stabilizer of the point
$\frakv$ of $\frak g^*$ under the coadjoint action; we note that
the stabilizer $G_p$ of $p$ is plainly a subgroup of $G_{\frakv}$.
Endow the Lie algebra $\frak g$ with a $G_{\qq}$-invariant inner
product, let $\frak g_{\frakv}\subseteq \frak g$ be the Lie
algebra of $G_{\frakv}$ and $\frak g_p \subseteq \frak g_{\frakv}$
that of $G_p$, let $\lambd$ be the orthogonal complement of $\frak
g_{\frakv}$ in $\frak g$, and let $\frak m$ be the orthogonal
complement of  $\frak g_{p}$ in $\frak g_{\frakv}$, so that
$$
\align \frak g &= \frak g_{\frakv} \oplus \lambd = \frak g_{p}
\oplus \frak m \oplus \lambd, \tag0.1.1
\\
\frak g^* &= \frak g_{\frakv}^* \oplus \lambd^* =\frak g^*_{p}
\oplus \frak m^* \oplus \lambd^* \tag0.1.2
\endalign
$$
are $G_{p}$-invariant orthogonal decompositions. Then
$$
\aligned \frak g p &= \frak g_{\frakv}p \oplus \lambd p =\frak m p
\oplus \lambd p \subseteq \roman T_p(\roman T^* \QQ),
\\
(\frak g p)^* &= (\frak g_{\frakv}p)^* \oplus (\lambd p)^* =(\frak
m p)^* \oplus (\lambd p)^*
\endaligned
\tag0.2
$$
are $G_p$-invariant decompositions, indeed, the projection from
$\frak g$ onto $\frak g p$ decomposes into the direct sum of the
surjection from $\frak g_{\frakv}$ onto $\frak g_{\frakv}p= \frak
m p $ with kernel $\frak g_p$ and an isomorphism from $\lambd$
onto $\lambd p$, whence the projection from $\frak g$ onto $\frak
g p$, restricted to $\frak m \oplus \lambd$, is an isomorphism
onto $\frak m p \oplus \lambd p = \frak g p$. In Section 4 below
we shall show that the cotangent bundle projection map induces a
$G_p$-equivariant isomorphism from $(\frak g \qq)^* \oplus (\frak
g_{\qq}p)^*$ onto $(\frak g p)^* =(\frak m p)^* \oplus (\lambd
p)^*$ which, in view of the decompositions (0.2), enables us to
decompose the infinitesimal $G$-slice $S_p$ given in the above
theorem as
$$
S_p \cong (\frak m p)^* \oplus (\lambd p)^*
 \oplus \sigma^*_{\qq}\oplus
\sigma_{\qq}
$$
and hence the tangent space $\roman T_p(\roman T^*\QQ)$ as
$$
 \frak g p \oplus S_p \cong \frak m p \oplus
\lambd p\oplus (\frak m p)^* \oplus (\lambd p)^*
 \oplus \sigma^*_{\qq}\oplus
\sigma_{\qq}.
$$
Now Lemma 4.8 below will say that the derivative $d\mu_p\colon
\roman T_p(\roman T^*\QQ) \to \frak g^*$ of the momentum mapping
$\mu$ at the point $p$ vanishes on the sum $\frak m p
\oplus\sigma^*_{\qq}\oplus \sigma_{\qq}$, comes down to the
canonical injection of $(\frak m p)^*\oplus (\lambd p)^*$ into
$\frak g^*$ on $(\frak m p)^*\oplus (\lambd p)^*$, the canonical
injection being induced by the projection from $\frak g$ onto
$\frak g p=\frak m p \oplus \lambd p\cong \frak g\big/ \frak g_p$
and, furthermore, amounts to the map
$$
\lambd \oplus \lambd ^* @>>> \lambd ^*,\quad (X,\alpha)
\longmapsto \alpha - \roman {ad}_X^*(\nu),\quad X \in \lambd,
$$
on $\lambd p\oplus (\lambd p)^*$. Let $K$ be the kernel of the
latter map. The vector space $(\lambd p)^*$ is canonically
isomorphic to the tangent space $\roman T_{\frakv}\Cal O_{\frakv}$
to the orbit $\Cal O_{\frakv}$ at $\frakv$ and the vector space
$\roman T_{\frakv}\Cal O_{\frakv}$, in turn, endowed with the
negative Kostant-Kirillov-Souriau symplectic form, is a symplectic
vector space. Now, under the projection from $\lambd p\oplus
(\lambd p)^*$ to the second summand, $K$ is mapped isomorphically
onto $(\lambd p)^*$ in such a way that the symplectic form on
$\roman T_p(\roman T^*\QQ)$, restricted to $K$, is identified with
the negative Kostant-Kirillov-Souriau symplectic form. Let
$V=K\oplus \sigma^*_{\qq}\oplus \sigma_{\qq}$, and endow $V$ with
the product symplectic structure, the structure on $
\sigma^*_{\qq}\oplus \sigma_{\qq}$ being the cotangent bundle
structure. The obvious $G_{\frakv}$-action on $K$ and the obvious
$G_{\qq}$-action on $\sigma^*_{\qq}\oplus \sigma_{\qq}$ induce a
$G_p$-action on $V$, necessarily Hamiltonian in the obvious
fashion. The Hamiltonian $G_p$-space $V$ is a {\it symplectic
slice\/} at the point $p$ of $\roman T^*\QQ$. Moreover, letting $W
= (\frak m p)^*$, we obtain the {\it Witt-Artin decomposition\/}
$$
\roman T_p(\roman T^* \QQ) = \frak g_{\nu} p \oplus \lambd p
\oplus V \oplus W \tag0.3
$$
at the point $p$ of $\roman T^* \QQ$. We shall show in Section 4
below how this decomposition can entirely be characterized in
terms of the $G$-action on the base $\QQ$ and of the geometry of
the coadjoint orbit generated by $\nu$.

In our approach, the group $G$ and manifold $\QQ$ may be infinite
dimensional. The only additional requirement is, then, that the
Riemannian metrics used below always exist as {\it strong\/}
metrics and that the implicit function theorem hold. This
situation actually arises in gauge theory.

The present paper was prompted by the recent postings of
\cite\perrosou\ and \cite\schmah.

\beginsection 1. Slices and fiber bundles

Let $\pi \colon E \to \QQ$ be a fiber bundle. For $\qq \in \QQ$,
let $F_{\qq}\subseteq E$ be the fiber $\pi^{-1}(\qq)$ over $\qq$.
The (total spaces of the) tangent bundles constitute a fiber
bundle
$$
\roman T\pi \colon \roman TE @>>> \roman T\QQ \tag1.1
$$
in an obvious fashion. Given a point $\qq$ of $\QQ$ and a vector
$v_{\qq}$ in the tangent space $\roman T_{\qq}\QQ$ at $\qq$, the
fiber $(\roman T\pi)^{-1}(v_{\qq})$ over $v_{\qq}$ is the bundle
$\roman AF_{\qq}$ of affine spaces over $F_{\qq}$ such that, for
$p\in F_{\qq}$, the fiber $\roman A_pF_{\qq}$ is the affine
subspace $w_p + \roman T_p F_{\qq}$ of $\roman T_p E$ where
$w_p\in \roman T_p E$ is a pre-image of $v_{\qq}$; the affine
subspace $w_p + \roman T_p F_{\qq}$ depends only on $v_{\qq}$ and
$p$ and not on the choice of $w_p$ since two such choices differ
by an element of $\roman T_p F_{\qq}$. Notice that when $v_{\qq}$
is the zero vector of $\roman T_{\qq}\QQ$ the fiber $(\roman
T\pi)^{-1}(v_{\qq})$ over $v_{\qq}$ is the ordinary tangent bundle
$\roman TF_{\qq}$ of $F_{\qq}$.
 In particular,
at $p\in E$ with $\pi(p) =\qq \in \QQ$, the tangent spaces
constitute the extension
$$
0 @>>> \roman T_pF_{\qq} @>>> \roman T_pE @>>> \roman T_{\qq}\QQ
@>>> 0 \tag1.2
$$
of vector spaces.

Let $G$ be a Lie group, write its Lie algebra as $\frak g$, and
suppose that $\pi$ is a $G$-fiber bundle. Thus $G$ acts on $E$ and
on $\QQ$ and the projection $\pi$ is a $G$-map. Given $(y,r)\in
G\times \QQ$ we will write the result of the action of $y$ on $r$
as $yr \in \QQ$ and, likewise, given $(y,z)\in G\times E$, we will
write the result of the action of $y$ on $z$ as $yz \in E$.

Let $p$ be a point of $E$, let $\qq=\pi(p)\in \QQ$, and suppose
that the stabilizer $G_{\qq}$ of $\qq$ is a compact Lie group.
Then the stabilizer $G_p$ of $p$ is a closed subgroup of $G_{\qq}$
and hence a compact Lie group as well. For example, when the
action is proper, the stabilizer $G_y$ is compact for any point
$y$ of $\QQ$. Our aim is to explore the relationship between
$G$-slices at $p$ and $\qq$.

To this end, we first endow $E$ and $\QQ$ with $G_{\qq}$-invariant
Riemannian metrics in such a way that, at each point $y$ of $E$,
the orthogonal complement of the tangent space  $\roman T_yF$ to
the fiber $F$  at $y$ is mapped under $\pi$ isometrically to the
tangent space $\roman T_{\pi(y)}\QQ$ to the base $\QQ$ at the
point $\pi(y)$ of $\QQ$; we refer to this situation by saying that
the projection $\pi$ is {\it compatible\/} with the Riemannian
structures. Such Riemannian metrics on $E$ and $\QQ$ can be
constructed as follows: Associated with $\pi$ is the standard
extension
$$
0 @>>> \verti @>>> \roman TE @>>> \pi^* \roman T\QQ @>>> 0 \tag
1.3
$$
of vector bundles on $E$ where $\verti$ denotes the (total space
of the) vertical subbundle, that is, the bundle of vectors tangent
to the fibers of $\pi$. At $p_r\in E$ with $\pi(p_r)=r \in \QQ$,
the extension (1.3) comes down to the extension
$$
0 @>>> \verti_{p_r} @>>> \roman T_{p_r}E @>>> \roman T_r\QQ @>>> 0
\tag1.4
$$
of vector spaces. This is just a rewrite of the extension (1.2),
with $r$ substituted for $\qq$ and $p_r$ for $p$. Recall that an
{\it Ehresmann connection\/} for $\pi$ is by definition a
splitting of the extension (1.3) of vector bundles on $E$. Such a
connection may be given by either a section $\pi^* \roman T\QQ \to
\roman TE$ of vector bundles over $E$, so that this section and
the inclusion of $\verti$ into $\roman TE$ induce an isomorphism
$$
\verti \oplus \pi^* \roman T\QQ  @>>> \roman TE \tag1.5
$$
of vector bundles over $E$, or by a surjection of vector bundles
$\roman TE @>>> \verti$ over $E$, referred to as a {\it connection
form\/}, so that this connection form and the projection from
$\roman TE$ onto $\pi^* \roman T\QQ$ induce an isomorphism
$$
\roman TE @>>> \verti \oplus \pi^* \roman T\QQ \tag1.6
$$
of vector bundles over $E$; as usual, the summand $\verti$ is then
referred to as the {\it vertical\/} part and the summand $\pi^*
\roman T\QQ$ as the {\it horizontal\/} part. When $M$ carries a
Riemannian metric, an Ehresmann connection for $\pi$ induces a
Riemannian metric on the total space $E$ of $\pi$ in such a way
that the decompositions (1.5) and (1.6) are orthogonal
decompositions at each point of $E$ and that the projection $\pi$
is compatible with the metrics in the sense explained earlier;
when the connection is, furthermore, compatible with the
$G_{\qq}$-module structures, the Riemannian metrics may be taken
to be $G_{\qq}$-invariant in such a way that the decompositions
(1.5) and (1.6) are orthogonal decompositions of
$G_{\qq}$-representations at each point of $E$.

Now, suppose that $\pi$ is endowed with a $G_{\qq}$-invariant
Ehresmann connection and suppose that $E$ and $\QQ$ are endowed
with corresponding $G_{\qq}$-invariant Riemannian metrics such
that $\pi$ is compatible with the metrics. Consider the induced
$G_{\qq}$-representation on $\roman T_{\qq} \QQ$ and let $S_{\qq}
= (\frak g \qq)^{\bot}$, the orthogonal complement of the tangent
space $\frak g \qq=\roman T_{\qq}(G\qq)$ to the $G$-orbit of $\qq$
in $\QQ$. This is the standard infinitesimal slice at $\qq$ for
the $G$-action on $\QQ$. Then a suitable $G_{\qq}$-invariant ball
$B_{\qq}\subseteq S_{\qq}$ containing the origin will be a local
slice, that is, the map
$$
G \times _{G_{\qq}} S_{\qq} @>>> \QQ, \quad (x,y) \mapsto x\cdot
\roman{exp}_{\qq}(y), \tag1.7
$$
restricted to $G \times _{G_{\qq}} B_{\qq}$, is a diffeomorphism
onto a $G$-invariant neighborhood of $\qq$ in $\QQ$ in such a way
that the zero section goes to the orbit $G\qq$; here
$\roman{exp}_{\qq}$ refers to the exponential mapping at the point
$\qq$ for the Riemannian metric on $\QQ$. We will refer to this
kind of situation as a {\it slice decomposition for the\/}
$G$-{\it action on\/} $\QQ$ {\it at the point\/} $\qq$. Likewise,
consider the induced $G_p$-representation on $\roman T_p E$ and
let $S_p = (\frak g p)^{\bot}$, the orthogonal complement in
$\roman T_p E$ of the tangent space $\frak g p=\roman T_p(Gp)$ to
the $G$-orbit of $p$ in $E$. This is the standard infinitesimal
slice at $p$ for the $G$-action on $E$, and the map
$$
G \times _{G_p} S_p @>>> E, \quad (x,y) \mapsto x
\roman{exp}_p(y), \tag1.8
$$
is a slice decomposition for the $G$-action on $E$ at the point
$p$; here $\roman{exp}_p$ refers to the exponential mapping at the
point $p$ for the Riemannian metric on $E$. By construction, the
projection from $\roman T_pE$ to $\roman T_{\qq}\QQ$ restricts to
a projection from $S_p$ onto $S_{\qq}$, the kernel $\kappa_p$
thereof inherits a $G_{\qq}$-module structure so that $\kappa_p$
is actually an orthogonal $G_{\qq}$-representation, and the
resulting map
$$
G_{\qq} \times _{G_p} \kappa_p  @>>> F_{\qq} \tag1.9
$$
is a slice decomposition for the $G_{\qq}$-action on the fiber
$F_{\qq}$ at the point $p$ of $F_{\qq}$; notice that, indeed,
$\kappa_p$ is the orthogonal complement $(\frak g_{\qq} p)^{\bot}$
of $\frak g_{\qq} p$ in $\roman T_pF_{\qq}$. Furthermore, the
obvious projection from $G \times _{G_p} S_p$ to $G \times
_{G_{\qq}} S_{\qq}$ is a fiber bundle having the space $G_{\qq}
\times _{G_p} \kappa_p$ as fiber at $\qq$, and the resulting
diagram
$$
\CD
G_{\qq} \times _{G_p} \kappa_p  @>>> F_{\qq}\\
@VVV @VVV
\\
G \times _{G_p} S_p  @>>> E\\
@VVV @VVV
\\
G \times _{G_{\qq}} S_{\qq} @>>> \QQ
\endCD
\tag1.10
$$
whose top and middle horizontal maps are the corresponding slice
decompositions at $p$ and whose bottom horizontal map is the
corresponding slice decomposition at $\qq$ is commutative, that
is, a morphism of fibre bundles. Now, consider the orthogonal
complement $\kappa_p^{\bot}$ of $\kappa_p$ in $S_p$ and the
orthogonal complement $(\frak g_{\qq} p)^{\bot}$ of $\frak
g_{\qq}p$ in $\frak gp$, viewed as a linear subspace of $\roman
T_pE$. In terms of these vector spaces, we have the decompositions
$$
\roman T_pF = \frak g_{\qq}p \oplus \kappa_p, \quad \roman T_pE =
\frak g_{\qq}p \oplus (\frak g_{\qq} p)^{\bot}\oplus \kappa_p
\oplus \kappa_p^{\bot}
$$
of $G_{\qq}$-representations and, as an extension of
$G_{\qq}$-representations, (1.2) may be written as
$$
0 @>>> \frak g_{\qq}p \oplus \kappa_p @>>>  \frak g_{\qq}p \oplus
(\frak g_{\qq} p)^{\bot}\oplus \kappa_p \oplus \kappa_p^{\bot}
@>>> \frak g\qq \oplus S_{\qq} @>>> 0 \tag1.11
$$
in such a way that the injection is the obvious one and  the
projection identifies $(\frak g_{\qq} p)^{\bot}$ with $\frak g\qq$
and $\kappa_p^{\bot}$ with $S_{\qq}$. A choice of suitable balls
in the infinitesimal slices then yields local slices.

\medskip\noindent{\bf 2. The tangent bundle}\smallskip\noindent
We now apply the previous discussion to the special case where
$\pi$ is the tangent bundle of $\QQ$ and where the total space is
endowed with the lifted $G$-action. Since below we shall explore
the tangent bundle of the total space of the tangent bundle,  for
clarity, we will write the tangent bundle of $\QQ$ as
$\tau_\QQ\colon \QQ^{\roman T}\to \QQ$, where $\QQ^{\roman T}$ is
the total space. As in the previous section, let $\qq$ be a point
of $\QQ$ with compact stabilizer $G_{\qq}$, and let $p$ be a
tangent vector at $\qq$, that is, a point of $\roman T_{\qq}\QQ$.
We will sometimes write the fiber $\roman T_{\qq}\QQ$ as
$F_{\qq}$. Furthermore, introduce $G_{\qq}$-invariant Riemannian
metrics in $\QQ$ and $\QQ^{\roman T}$ in such a way that the
projection $\tau_\QQ$ is compatible with the metrics in a sense
explained in the previous section. To this end, endow $\QQ$ with a
$G_{\qq}$-invariant Riemannian metric. The corresponding
Levi-Civita connection on $\QQ$ is $G_{\qq}$-invariant and induces
a $G_{\qq}$-invariant Ehresmann connection on the tangent bundle
$\tau_\QQ$  in a canonical fashion. This connection, in turn, in
particular splits the corresponding exact sequence of orthogonal
$G_{\qq}$-representation spaces of the kind (1.4), that is,
induces a $G_{\qq}$-equivariant isomorphism
$$
\varphi^{\roman T}\colon\roman T_p(\roman T_{\qq}\QQ) \oplus
\roman T_{\qq}\QQ @>>> \roman T_{p}(\QQ^{\roman T}) \tag2.1
$$
of orthogonal $G_{\qq}$-representations.

Consider the infinitesimal slices $S_{\qq} = (\frak g\qq)^{\bot}
\subseteq \roman T_{\qq}\QQ$, $S_p = (\frak gp)^{\bot} \subseteq
\roman T_p\QQ^{\roman T}$, and $\kappa_p = (\frak g_{\qq}p)^{\bot}
\subseteq \roman T_pF_{\qq}$ as well as the slice decompositions
$$
\align G &\times _{G_{\qq}} S_{\qq} @>>> \QQ , \tag2.2
\\
G &\times_{G_p} S_p  @>>> \QQ^{\roman T}, \tag2.3
\\
G_{\qq} &\times_{G_p} \kappa_p  @>>> F_{\qq} \tag2.4
\endalign
$$
introduced in the previous section. Since, with reference to the
decomposition $\roman T_{\qq}M = \frak g \qq \oplus S_{\qq}$, the
$G_p$-action on $\frak g \qq$ is trivial, the tangent space $\frak
g_{\qq}p = \roman T_p(G_{\qq}p) \subseteq \roman T_{\qq}M$ to the
$G_{\qq}$-orbit of $p$ in $\roman T_{\qq}M$ is actually  a linear
subspace of $\roman T_pS_{\qq}\cong S_{\qq}$; let $\sigma_{\qq}
\subseteq S_{\qq}$ be the orthogonal complement of $\frak
g_{\qq}p$ in the vector space $S_{\qq}$. Then $S_{\qq} =\frak
g_{\qq}p \oplus \sigma_{\qq}$ and hence
$$
\roman T_{\qq} \QQ = \frak g \qq \oplus S_{\qq} = \frak g \qq
\oplus \frak g_{\qq} p \oplus \sigma_{\qq} \tag2.5
$$
are $G_p$-invariant orthogonal decompositions. Accordingly, the
isomorphic image $\varphi^{\roman T}(\roman T_{\qq} \QQ)$
decomposes as
$$
\varphi^{\roman T}(\roman T_{\qq} \QQ) = (\frak g \qq)^h \oplus
(\frak g_{\qq} p)^h \oplus \sigma^h_{\qq} \tag2.6
$$
in such a way that the summands $(\frak g \qq)^h$, $(\frak g_{\qq}
p)^h$ and $\sigma^h_{\qq}$ correspond to, respectively, $\frak g
\qq$, $\frak g_{\qq} p$, and $\sigma_{\qq}$,  whence the notation
$(\frak g \qq)^h$, $(\frak g_{\qq} p)^h$ and $\sigma^h_{\qq}$;
here the superscript $h$ is intended to indicate that, in the
decomposition (2.9.1) below, $(\frak g p)^h$ and $\sigma^h_{\qq}$
will be {\it horizontal\/} constituents with respect to the
Ehresmann connection.

Under the canonical isomorphism of vector spaces, even orthogonal
$G_{\qq}$-re\-pre\-sen\-ta\-tions, between $\roman T_{\qq}\QQ$ and
 $\roman T_p(F_{\qq})=\roman T_p(\roman
T_{\qq}\QQ)$, the orthogonal decomposition (2.5) corresponds to
the decomposition
$$
\roman T_p(F_{\qq})= (\frak g \qq)^v \oplus (\frak g_{\qq} p)^v
\oplus \sigma^v_{\qq} \tag2.7
$$
in such a way that the summands $(\frak g \qq)^v$, $(\frak g_{\qq}
p)^v$ and $\sigma^v_{\qq}$ correspond to, respectively, $\frak g
\qq$,  $\frak g_{\qq} p$, and $\sigma_{\qq}$, whence the notation
$(\frak g \qq)^v$, $(\frak g_{\qq} p)^v$ and $\sigma^v_{\qq}$;
here the superscript $v$ is intended to indicate that, in the
decomposition (2.9.1) below, $(\frak g p)^v$ and $\sigma^v_{\qq}$
will be {\it vertical\/} constituents with respect to the
Ehresmann connection. Then
$$
\kappa_p = (\frak g \qq)^v \oplus \sigma^v_{\qq}, \quad
\kappa^{\bot}_p = (\frak g_{\qq}p)^h \oplus \sigma^h_{\qq},
 \tag2.8
$$
and the construction in the previous section entails the
following.

\proclaim{Theorem 2.9} The infinitesimal slice $S_p = (\frak g
p)^{\bot} \subseteq \roman T_p(\QQ^{\roman T})$ has the
$G_p$-invariant orthogonal decomposition
$$
S_p = (\frak g \qq)^v \oplus \sigma^v_{\qq} \oplus (\frak
g_{\qq}p)^h \oplus \sigma^h_{\qq}. \tag2.9.1
$$
Furthermore, the projection from $S_p$ to $S_{\qq}$, restricted to
$\sigma^v_{\qq} \oplus\sigma^h_{\qq}$, coincides with the obvious
projection from $\sigma^v_{\qq} \oplus\sigma^h_{\qq}$ to
$\sigma^h_{\qq}$ followed by the canonical isomorphism from
$\sigma^h_{\qq}$ onto the summand $\sigma_{\qq}$ of $S_{\qq}$ and
hence amounts to the tangent bundle projection map of
$\sigma_{\qq}$.
\endproclaim

\demo{Proof} Indeed, $S_p = \kappa_p \oplus \kappa_p^{\bot}$ and,
as explained in the previous section, $\frak g p$ decomposes as
the direct sum $\frak g_{\qq} p \oplus (\frak g_{\qq} p)^{\bot}$
of $\frak g_{\qq} p$ with its orthogonal complement $(\frak
g_{\qq} p)^{\bot}$ in $\frak g p$ in such a way that, under the
tangent bundle projection map, $(\frak g_{\qq} p)^{\bot}$ is
identified with $\frak g {\qq}$. This implies the assertion. \qed
\enddemo

Turned the other way round, the theorem says the following.

\proclaim{Corollary 2.10} At the point $p$ of $\roman
T_{\qq}(\QQ)$, as an orthogonal $G_p$-representation, a suitable
infinitesimal slice for the lifted $G$-action on the total space
$\QQ^{\roman T}$ of the tangent bundle of $\QQ$ is isomorphic to
the orthogonal $G_p$-representation
$$
\frak g \qq\oplus \sigma_{\qq} \oplus \frak g_{\qq}p \oplus
\sigma_{\qq}.
$$
\endproclaim

This description of the infinitesimal $G$-slice at the point $p$
of the total space $\QQ^{\roman T}$ of the tangent bundle of $\QQ$
plainly involves  only data phrased in terms of the base $\QQ$.

The above slice decomposition admits an interpretation in terms of
the induced action of the tangent group. We now explain this,
since it will help understand the cotangent bundle situation in
the next section. Thus, we write the total space of the tangent
bundle of $G$ as $G^{\roman T}$ and note first that $G^{\roman T}$
inherits a Lie group structure in the following fashion: Let
$\frak g_a$ be the Lie algebra $\frak g$, viewed as an abelian Lie
algebra, i.~e. real vector space, endowed with the
$G$-representation coming from the adjoint action. Then the
semi-direct product Lie group $G \ltimes \frak g_a$ is defined,
and left translation identifies $G \ltimes \frak g_a$ with
$G^{\roman T}$ and hence turns $G^{\roman T}$ into a Lie group in
such a way that the tangent bundle projection map is a surjective
homomorphism. Explicitly, given $V_1, V_2 \in \frak g$ and
$y_1,y_2 \in G$, write $y_1V_1 \in \roman T_{y_1}G$ and $y_2V_2
\in \roman T_{y_2}G$ for the results of left translation with
$y_1$ and $y_2$, respectively; then
$$
(y_1 V_1)(y_2 V_2) = y_1y_2\left(\roman{Ad}(y_2)^{-1}V_1 +
V_2\right) \in \roman T_{y_1y_2}G. \tag2.11
$$
Given $r \in \QQ$ and $(y,X_r)\in G \times \roman T_r\QQ$, with
reference to the lifted $G$-action $ G\times \QQ^{\roman T} \to
\QQ^{\roman T}$ on $\QQ^{\roman T}$, we will write the result of
the action of $y$ on $X_r$ as $ yX_r \in \roman T_{yr}\QQ$. The
lifted action, in turn, extends to a $G^{\roman T}$-action
$$
G^{\roman T} \times \QQ^{\roman T} @>>> \QQ^{\roman T} \tag2.12
$$
on  $\QQ^{\roman T}$ in an obvious manner: Given $y \in G$, $V\in
\frak g$, $r\in \QQ$, and $X_r \in \roman T_r\QQ \subseteq
\QQ^{\roman T}$, the result of the action of $yV \in G^{\roman T}$
upon $X_r$ is given by the expression
$$
yV(X_r) = y X_r + y V_r. \tag2.13
$$
Here $V_r \in \roman T_r\QQ$ refers to the value at $r \in \QQ$ of
the fundamental vector field $V_\QQ$ on $\QQ$ induced by $V\in
\frak g$, and $y X_r + y V_r$ arises from the pair $(y, X_r + V_r)
\in G \times \roman T_r \QQ$ by an application of the lifted
$G$-action on $\QQ^{\roman T}$. It is straightforward to check
that (2.13) yields indeed an action of $G^{\roman T}$ on
$\QQ^{\roman T}$. In particular, the restriction of the action to
the copy of $\frak g_a$ maps each fiber $\roman T_r\QQ\subseteq
\QQ^{\roman T}$ to itself where $r \in \QQ$, and the resulting
action
$$
\frak g_a \times \roman T_r\QQ\subseteq \QQ^{\roman T} @>>> \roman
T_r\QQ\subseteq \QQ^{\roman T} \tag2.14
$$
of $\frak g_a$ on $\roman T_r\QQ$ is given by the assignment to
$(V, X_r) \in \frak g_a \times \roman T_r\QQ\subseteq \QQ^{\roman
T}$ of $X_r + V_r\in \roman T_r\QQ$ and is hence given by affine
transformations, i.~e. generalized translations. Consequently, for
each $r \in \QQ$, the induced infinitesimal action
$$
\frak g_a @>>> \roman{Vect}(\roman T_r\QQ) \tag2.15
$$
sends each $W\in \frak g_a$ to the constant vector field $W_r$
which assigns to any point $X_r \in \roman T_r\QQ$ the vector $W_r
\in \roman T_{X_r}(\roman T_r\QQ)$, the latter vector space being
canonically identified with $\roman T_r\QQ$, viewed as a vector
space.

Consider the slice decomposition (2.2), and let
$\tau_{G_{\qq}}\colon G^{\roman T}_{\qq} \to G_{\qq}$ and
$\tau_{S_{\qq}}\colon S^{\roman T}_{\qq} \to S_{\qq}$ be the
tangent bundles of $G_{\qq}$ and $S_{\qq}$, respectively.
Similarly as above, left translation identifies the semi-direct
product $G_{\qq} \ltimes (\frak g_{\qq})_a$ with $G_{\qq}^{\roman
T}$. The resulting slice decomposition
$$
G^{\roman T} \times _{G^{\roman T}_{\qq}} S^{\roman T}_{\qq} @>>>
\QQ^{\roman T}, \quad (x,y) \mapsto x (\roman
T\roman{exp}_{\qq})(y), \tag2.16
$$
for the $G^{\roman T}$-action on $\QQ^{\roman T}=\roman T\QQ$ at
the point $0_{\qq}$ of $\roman T\QQ$, that is, at the origin of
the tangent space $\roman T_{\qq}\QQ$ of $\qq$ at $\QQ$, viewed as
a subspace of $\QQ^{\roman T}=\roman T\QQ$, exhibits the total
space of the tangent bundle of $G \times _{G_{\qq}} S_{\qq}$ as
the quotient space $G^{\roman T} \times _{G^{\roman T}_{\qq}}
S^{\roman T}_{\qq}$ in such a way that this slice decomposition
amounts to the tangent map of the slice decomposition $G \times
_{G_{\qq}} S_{\qq} @>>> \QQ$ for the $G$-action on the base $\QQ$
at the point $\qq$, cf. (2.2).

\medskip\noindent{\bf 3. The cotangent bundle}\smallskip\noindent
We now explain the necessary modifications for the case where
$\pi$ is the cotangent bundle $\tau^*_\QQ\colon\roman T^*\QQ \to
\QQ$ of $\QQ$ and where the total space is endowed with the lifted
action of the Lie group $G$. Without further hypotheses, we cannot
simply identify the tangent bundle with the cotangent bundle in a
$G$-equivariant fashion, and we must be a bit more circumspect.

Below it will be convenient to downplay the structure of $\roman
T^*\QQ$ as the total space of a vector bundle and to concentrate
merely on the structure as a symplectic manifold. We will then
denote the total space $\roman T^*\QQ$ by $\QQ^{\roman T^*}$. Let
$\qq$ be a point of $\QQ$ with compact stabilizer $G_{\qq}$, and
let $p$ be a cotangent vector at $\qq$, that is, a point of
$\QQ^{\roman T^*}$ with $\tau^*_\QQ(p) = \qq$. We will sometimes
write the fiber $\roman T^*_{\qq}\QQ$ as $F_{\qq}$.

Similarly as in the previous section, the Levi-Civita connection
on $\QQ$ induced by the chosen  $G_{\qq}$-invariant Riemannian
metric on $\QQ$ is $G_{\qq}$-invariant and induces a
$G_{\qq}$-invariant Ehresmann connection on the cotangent bundle
$\tau^*_\QQ$  in a canonical fashion. This connection, in turn, in
particular splits the corresponding exact sequence of orthogonal
$G_{\qq}$-representation spaces of the kind (1.4), that is,
induces a $G_{\qq}$-equivariant isomorphism
$$
\varphi^{\roman T^*}\colon\roman T_p(\roman T^*_{\qq}\QQ) \oplus
\roman T_{\qq}\QQ @>>> \roman T_{p}(\QQ^{\roman T^*}) \tag3.1
$$
of orthogonal $G_{\qq}$-representations.  Consider the
infinitesimal slices $S_{\qq} = (\frak g\qq)^{\bot} \subseteq
\roman T_{\qq}\QQ$, $S_p = (\frak gp)^{\bot} \subseteq \roman
T_p\QQ^{\roman T^*}$, and $\kappa_p = (\frak g_{\qq}p)^{\bot}
\subseteq \roman T_pF_{\qq}$ as well as the slice decompositions
$$
\align G &\times _{G_{\qq}} S_{\qq} @>>> \QQ , \tag3.2
\\
G &\times_{G_p} S_p  @>>> \QQ^{\roman T^*}, \tag3.3
\\
G_{\qq} &\times_{G_p} \kappa_p  @>>> F_{\qq}, \tag3.4
\endalign
$$
introduced in Section 1 above. As in the previous section, our
present aim is now to decompose $S_p$ and $S_{\qq}$ further in
such a way that a summand of $S_p$ may be identified with the
total space of the tangent bundle of a summand of $S_{\qq}$.

To this end we note first that, even though there is no obvious
way to identify the tangent and cotangent bundles of $\QQ$ in a
$G$-equivariant manner, the chosen $G_{\qq}$-invariant Riemannian
metric on $\QQ$ induces a $G_{\qq}$-equivariant isomorphism
between the vector spaces $\roman T_{\qq}\QQ$ and $\roman
T_{\qq}^*\QQ$, and the $G_{\qq}$-invariant inner product on
$\roman T_{\qq}\QQ$ induces a $G_{\qq}$-invariant inner product on
$\roman T_{\qq}^*\QQ$. Thus the orthogonal decomposition $\roman
T_{\qq} \QQ = \frak g \qq \oplus S_{\qq}$ induces the orthogonal
decomposition $\roman T^*_{\qq} \QQ = (\frak g \qq)^* \oplus
S^*_{\qq}$.

Consider the tangent space $\frak g_{\qq}p = \roman T_p(G_{\qq}p)
\subseteq \roman T^*_{\qq}M$ to the $G_{\qq}$-orbit of $p$ in
$\roman T^*_{\qq}M$. Since, with reference to the decomposition
$\roman T^*_{\qq}M = (\frak g \qq)^* \oplus S^*_{\qq}$, the
$G_p$-action on $(\frak g \qq)^*$ is trivial, the tangent space
$\frak g_{\qq}p = \roman T_p(G_{\qq}p) \subseteq \roman
T^*_{\qq}M$ to the $G_{\qq}$-orbit of $p$ in $\roman T^*_{\qq}M$
is actually  a linear subspace of $\roman T_pS^*_{\qq}\cong
S^*_{\qq}$; let $\sigma_{\qq} \subseteq S_{\qq}$ be the unique
linear subspace such that $S^*_{\qq} =\frak g_{\qq}p \oplus
\sigma^*_{\qq}$ is a $G_p$-invariant orthogonal decomposition. Let
$(\frak g_{\qq}p)^*$ be the unique subspace of $S_{\qq}$ such that
its dual $(\frak g_{\qq}p)^{**}$ coincides with $\frak
g_{\qq}p\subseteq S^*_{\qq}$. Then $S_{\qq} =(\frak g_{\qq}p)^*
\oplus \sigma_{\qq}$ is a $G_p$-invariant orthogonal decomposition
as well, and
$$
\align \roman T^*_{\qq} \QQ &= (\frak g \qq)^* \oplus S^*_{\qq} =
(\frak g \qq)^* \oplus \frak g_{\qq} p \oplus \sigma^*_{\qq}
\tag3.5.1
\\
\roman T_{\qq} \QQ &= \frak g \qq \oplus S_{\qq} = \frak g \qq
\oplus (\frak g_{\qq} p)^* \oplus \sigma_{\qq} \tag3.5.2
\endalign
$$
are $G_p$-invariant orthogonal decompositions. Accordingly, the
isomorphic image $\varphi^{\roman T^*}(\roman T_{\qq} \QQ)$
decomposes as
$$
\varphi^{\roman T^*}(\roman T_{\qq} \QQ) = \frak g \qq \oplus
(\frak g_{\qq} p)^* \oplus \sigma_{\qq} \tag3.6
$$
where the notation $\frak g \qq$, $(\frak g_{\qq} p)^*$, and
$\sigma_{\qq}$ is slightly abused.

Under the canonical isomorphism of vector spaces, even orthogonal
$G_{\qq}$-re\-pre\-sen\-ta\-tions, between $\roman T^*_{\qq}\QQ$
and
 $\roman T_p(F_{\qq})=\roman T_p(\roman
T^*_{\qq}\QQ)$, the orthogonal decomposition (3.5.2) passes to the
decomposition
$$
\roman T_p(F_{\qq})=(\frak g \qq)^* \oplus \frak g_{\qq} p \oplus
\sigma^*_{\qq} \tag3.7
$$
where the notation $(\frak g \qq)^*$, $\frak g_{\qq} p$ and
$\sigma^*_{\qq}$ is slightly abused. Then
$$
\kappa_p = (\frak g \qq)^* \oplus \sigma^*_{\qq},
 \quad \kappa^{\bot}_p = (\frak g_{\qq}p)^* \oplus \sigma_{\qq},
 \tag3.8
$$
and the construction in Section 1 entails the following.

\proclaim{Theorem 3.9} The infinitesimal slice $S_p =(\frak g
p)^{\bot}\subseteq \roman T_p(\QQ^{\roman T^*})$ has the
$G_p$-invariant orthogonal decomposition
$$
S_p = (\frak g \qq)^* \oplus \sigma^*_{\qq}\oplus (\frak
g_{\qq}p)^* \oplus \sigma_{\qq}. \tag3.9.1
$$
Furthermore, the projection from $S_p$ to $S_{\qq}$, restricted to
$\sigma^*_{\qq}\oplus\sigma_{\qq}$, coincides with the obvious
projection from $\sigma^*_{\qq} \oplus\sigma_{\qq}$ to
$\sigma_{\qq}$ and hence amounts to the cotangent bundle
projection map of $\sigma_{\qq}$.
\endproclaim

We note that, in the statement of Theorem 3.9, there is no need to
distinguish in notation between horizontal and vertical
constituents since by construction $(\frak g \qq)^*$ and
$\sigma^*_{\qq}$ are vertical and $(\frak g_{\qq}p)^*$ and
$\sigma_{\qq}$ are horizontal with respect to the Ehresmann
connection.

\demo{Proof} The argument is formally the same as that for the
proof of  Theorem 2.9: The slice $S_p$ decomposes as $S_p=
\kappa_p \oplus \kappa_p^{\bot}$ and, as explained in Section 1,
$\frak g p$ decomposes as the direct sum $\frak g_{\qq} p \oplus
(\frak g_{\qq} p)^{\bot}$ of $\frak g_{\qq} p$ with its orthogonal
complement $(\frak g_{\qq} p)^{\bot}$ in $\frak g p$ in such a way
that, under the cotangent bundle projection map, $(\frak g_{\qq}
p)^{\bot}$ is identified with $\frak g {\qq}$. This implies the
assertion. \qed
\enddemo

Turned the other way round, Theorem 3.9 entails the theorem
spelled out in the introduction; the theorem in the introduction
is the exact analogue of Corollary 2.10 above.

\smallskip
\noindent {\smc Remark 3.10.} Let $X \in \roman T_{\qq}\QQ$ and
let $p \in \roman T^*_{\qq}\QQ$ be the image of $X$ under the
adjoint of the chosen $G_{\qq}$-invariant Riemannian metric. This
metric induces an isomorphism of $G_{\qq}$-representations from
the infinitesimal $G$-slice
$$
S_X= (\frak g \qq)^v \oplus \sigma^v_{\qq} \oplus (\frak
g_{\qq}X)^h \oplus \sigma^h_{\qq}
$$
given as (2.9.1) above onto the infinitesimal $G$-slice
$$
S_p = (\frak g \qq)^* \oplus \sigma^*_{\qq}\oplus (\frak
g_{\qq}p)^* \oplus \sigma_{\qq},
$$
given as (3.9.1) above. By construction, this isomorphism
decomposes into the direct sum of the induced isomorphisms
$$
(\frak g \qq)^v \to (\frak g \qq)^*, \ \sigma^v_{\qq} \to
\sigma^*_{\qq}, \ (\frak g_{\qq}X)^h \to (\frak g_{\qq}p)^*
$$
with the identity morphism  $\sigma^h_{\qq}\to \sigma_{\qq}$.
Moreover, the Riemannian metric identifies the stabilizers $G_X$
and $G_p$. In view of the slice decompositions (2.3) and (3.3),
the induced map from $ G \times _{G_X} S_X$ to $G \times _{G_p}
S_p$ therefore yields a $G$-equivariant diffeomorphism from a
neighborhood of the $G$-orbit of the point $X$ of $M^{\roman T}$
onto a neighborhood of the $G$-orbit of the point $p$ of
$M^{\roman T^*}$. From this observation it is straightforward to
concoct a proof of the familiar fact that, provided the group $G$
is finite dimensional, given a finite dimensional paracompact
$G$-manifold $Q$ having the property that the stabilizer of each
point is compact, the tangent and cotangent bundles of $Q$ are
$G$-equivariantly isomorphic as vector bundles on $Q$ or,
equivalently, the manifold $Q$ admits a $G$-invariant Riemannian
metric.

\beginsection 4. Symplectic slices and Witt-Artin decomposition

The purpose of this section is to give an interpretation of the
decomposition (3.9.1) in terms of a symplectic slice. We will also
construct a Witt-Artin decomposition, cf. \cite\ortrathr. Thus, as
before, let $G$ be a Lie group, $\QQ$ a $G$-manifold, and consider
the lifted $G$-action on $\QQ^{\roman T^*}\, (=\roman T^*\QQ)$.
Let $p$ be a point of $\QQ^{\roman T^*}$, and let $\qq = \pi(p)$
be the image of $p$ in $\QQ$ under the cotangent bundle projection
map $\pi = \tau_\QQ^*$. The notation $\mu$, $\frakv$, $\frak
g_{\frakv}$, $\lambd$ and $\frak m$ used below is the same as that
in the introduction.

Denote the infinitesimal coadjoint $\frak g$-action on $\frak g^*$
by $[\,\cdot , \cdot ] \colon \frak g \times \frak g^* \to \frak
g^*$, so that
$$
[Y,\psi](Z) = \psi([Z,Y]). \tag4.1
$$
It is also common in the literature to write
$(\roman{ad}_Y^*(\psi))(Z) = \psi([Y,Z])$. With the notation
$[\,\cdot , \cdot ]$,  the tangent map
$$
\frak g \to \roman T_{\frakv}\Cal O_{\frakv}\subseteq \frak g^*
\tag4.2
$$
of the canonical projection from $G$ onto $\Cal O_{\frakv} =
G\frakv$ is given by the assignment to $X\in \frak g$ of
$[X,\frakv]\in \frak g^*$. This tangent map passes to an
isomorphism from $\lambd$ onto $\roman T_{\frakv}\Cal O_{\frakv}$
and hence induces an isomorphism from $\lambd p$ onto $\roman
T_{\frakv}\Cal O_{\frakv}=(\lambd p)^*$---this isomorphism is
precisely the inverse of the adjoint of the symplectic structure
on $\roman T_{\frakv}\Cal O_{\frakv}$ corresponding to the
Kostant-Kirillov-Souriau symplectic structure on $\Cal
O_{\frakv}$. Consequently the tangent map of the canonical
projection from $G$ onto $\Cal O_{\frakv}$ induces
 an injection
$$
\lambd p @>>> \frak g^*. \tag4.3
$$
Furthermore, the projection
from $\frak g$ to $\frak g p$ induces an injection
$$
(\frak m p)^* \oplus (\lambd p)^* =(\frak g_{\frakv}p)^* \oplus
(\lambd p)^* =(\frak g p)^* @>>> \frak g^* . \tag4.4
$$
In view of Theorem 3.9, at the point $p$ of $\QQ^{\roman T^*}$,
the tangent space  $\roman T_p(\QQ^{\roman T^*})$ ($\cong\roman
T_p(G\times_{G_p} S_p)$) decomposes as
$$
\frak g p \oplus S_p = \frak g p \oplus  (\frak g \qq)^* \oplus
\sigma^*_{\qq}\oplus (\frak g_{\qq}p)^* \oplus \sigma_{\qq}.
$$
Consider the subspace $\frak g_{\qq}p$ of the tangent space $\frak
gp$ at $p$ to the $G$-orbit $Gp$ in $\QQ^{\roman T^*}$. With
respect to the $G_{\qq}$-invariant inner product on $\roman
T_p(\QQ^{\roman T^*})$,
 $\frak gp$ decomposes as
$$
\frak gp =\frak g_{\qq}p \oplus (\frak g_{\qq}p)^{\bot} \tag4.5
$$
in such a way that the projection from $\frak gp$ to $\frak g\qq$
(induced by the cotangent bundle projection map) restricts to an
isomorphism from $(\frak g_{\qq}p)^{\bot}$ onto $\frak g\qq$.
 This isomorphism, in turn, induces an isomorphism
$$
(\frak g \qq)^* \oplus (\frak g_{\qq}p)^* @>>> (\frak g p)^*
=(\frak m p)^* \oplus (\lambd p)^* \tag4.6
$$
which, in view of the decompositions (4.2), enables us to
decompose the infinitesimal $G$-slice $S_p$ as
$$
S_p \cong (\frak m p)^* \oplus (\lambd p)^*
 \oplus \sigma^*_{\qq}\oplus
\sigma_{\qq} \tag4.7.1
$$
and hence the tangent space $\roman T_p(\QQ^{\roman T^*})$ as
$$
 \frak g p \oplus S_p \cong \frak m p \oplus
\lambd p\oplus (\frak m p)^* \oplus (\lambd p)^*
 \oplus \sigma^*_{\qq}\oplus
\sigma_{\qq}. \tag4.7.2
$$

\proclaim{Lemma 4.8} In terms of the decomposition given as the
right-hand side of {\rm (4.7.2)}, at the point $p$, the derivative
$$
d\mu_p\colon \frak m p \oplus \lambd p\oplus (\frak m p)^* \oplus
(\lambd p)^*
 \oplus \sigma^*_{\qq}\oplus
\sigma_{\qq} @>>> \frak g^*
$$
of the momentum mapping $\mu$ vanishes on the sum $\frak m p
\oplus\sigma^*_{\qq}\oplus \sigma_{\qq}$, comes down to the
injection {\rm (4.3)} on the summand $\lambd p$, and to the
injection {\rm (4.4)} on $(\frak m p)^* \oplus (\lambd p)^* =
(\frak gp)^*$.
\endproclaim

To prepare for the proof, we note first that, since any $G$-slice
decomposition is a diffeomorphism onto a full neighborhood of a
$G$-orbit, the tangent map $\roman T(G\times_{G_{\qq}} S_{\qq})
\to \roman TM$ of the slice decomposition $G \times _{G_{\qq}}
S_{\qq} @>>> \QQ$ on the base, cf. (3.2), induces a smooth map
$$
\roman T^*(G\times_{G_{\qq}} S_{\qq}) @>>>  \roman T^*M =
M^{\roman T^*}
$$
over the slice decomposition on the base which is compatible with
the cotangent bundle structures. In particular, the induced map is
compatible with the tautological 1-forms and hence with the
symplectic structures and momentum mappings. We may therefore
replace the space $M^{\roman T^*}$ with the total space $\roman
T^*(G\times_{G_{\qq}} S_{\qq})$ of the cotangent bundle on
$G\times_{G_{\qq}} S_{\qq}$.

Write the total space of the cotangent bundle of the infinitesimal
slice $S_{\qq}$ as $S^{\roman T^*}_{\qq}$ and, likewise, write the
total space  of the cotangent bundle of $G$ as $G^{\roman T^*}$.
It is well known that, as a Hamiltonian $G$-space, a space of the
kind $\roman T^*(G\times_{G_{\qq}} S_{\qq})$ is in a canonical way
isomorphic to the Hamiltonian $G$-space $\left(G^{\roman
T^*}\times S^{\roman T^*}_{\qq}\right)_0$ which arises by
$G_{\qq}$-reduction at zero momentum, applied to the product
$G^{\roman T^*}\times S^{\roman T^*}_{\qq}$, endowed with the
product cotangent bundle symplectic structure, where the
$G_{\qq}$-action on $G^{\roman T^*}$ is induced from the lift of
the $G$-action on itself by {\it right\/} translation and where
the $G_{\qq}$-action on $S^{\roman T^*}_{\qq}$ is the cotangent
bundle lift of the $G_{\qq}$-action on $S_{\qq}$. The requisite
momentum mapping is, then, the sum of the momentum mappings of the
factors, where \lq\lq sum\rq\rq\ is interpreted in the obvious
fashion. In this way, the $G$-action on $G^{\roman T^*}$ which is
the lift of the $G$-action on itself by {\it left\/} translation
descends to a Hamiltonian $G$-action on the left of
$\left(G^{\roman T^*}\times S^{\roman T^*}_{\qq}\right)_0$ with
momentum mapping
$$
\mu_{[0]} \colon \left(G^{\roman T^*}\times S^{\roman
T^*}_{\qq}\right)_0 @>>> \frak g^*\tag4.8.1
$$
in such a way that the $G$-action and momentum mapping correspond
exactly to the corresponding cotangent bundle structure, with
reference to these structures on $\roman T^*(G\times_{G_{\qq}}
S_{\qq})$. This kind of construction underlies the local normal
form for the momentum mapping developed in \cite\guisteft\  and
\cite\marleone.

The above consideration reduces the proof of Lemma 4.8 to the
special case where the Hamiltonian $G$-space space $\QQ^{\roman
T^*}$ is of the kind $\left(G^{\roman T^*}\times S^{\roman
T^*}_{\qq}\right)_0$. The description of $\roman T^*(G\times
_{G_{\qq}}S_{\qq})$ as the reduced space $\left(G^{\roman
T^*}\times S^{\roman T^*}_{\qq}\right)_0$ is the exact dual of the
construction $G^{\roman T} \times _{G^{\roman T}_{\qq}} S^{\roman
T}_{\qq}$ in Section 2 above and the induced map
$$
\left(G^{\roman T^*}\times S^{\roman T^*}_{\qq}\right)_0 @>>>
\QQ^{\roman T^*} \tag4.8.2
$$
is the exact dual of (2.16) above.

\demo{Proof of Lemma 4.8} Under the slice decomposition $G \times
_{G_{\qq}} S_{\qq} @>>> \QQ$ for the $G$-action on the base $\QQ$
at the point $\qq$, cf. (3.2), the point $b$ of $\QQ$ corresponds
to the point $[e,0] \in G \times _{G_{\qq}}S_{\qq}$; here $[e,0]$
denotes the class in $G \times _{G_{\qq}}S_{\qq}$ of the point
$(e,0) \in G \times S_{\qq}$ where $e \in G$ is the neutral
element. Accordingly, the point $p$ now corresponds to a point of
the kind
$$
[e,\nu,0,p_{\qq}] \in \left(G^{\roman T^*}\times S^{\roman
T^*}_{\qq}\right)_0,
$$
that is, to the class in $\left(G^{\roman T^*}\times S^{\roman
T^*}_{\qq}\right)_0$ of a point of the kind
$$
(e,\nu,0,p_{\qq}) \in G^{\roman T^*}\times S^{\roman T^*}_{\qq}
$$
in the zero locus of the corresponding $G_{\qq}$-momentum mapping.
Recall that the familiar $G$-equivariant cotangent bundle momentum
mapping ${ \mu\colon G^{\roman T^*} @>>> \frak g^* }$ is  given by
the formula
$$
\mu(\alpha_x) =  \alpha_x \circ (R_x)_*,\quad x
\in G, \, \alpha_x \in \roman T^*_{x}G, 
$$
where, for $x \in G$, $(R_x)_*\colon \frak g = \roman T_eG
\to\roman T_{x}G$ refers to right translation by $x\in G$. In
particular, since the $G$-momentum mapping $\mu_{[0]}$ is induced
by $\mu$, necessarily $\mu_{[0]}[e,\nu,0,p_{\qq}] = \nu \in \frak
g^*$. The derivative
$$
d \mu_{\alpha}\colon \roman T_{e}G \times \roman T_\alpha \frak
g^* @>>> \frak g^*
$$
of the momentum mapping $\mu$ at the point $\alpha=\alpha_e\in
\roman T^*_eG =\frak g^*$ has the form
$$
d \mu_{\alpha}(X,\varphi)= \varphi + [X,\alpha ]. \tag4.8.3
$$
This implies the claim of Lemma 4.8. Indeed, in the special case
where $G_b$ is trivial, the decomposition (4.7.2) has the form
$$
\roman T_p(G\times S_p) \cong
 \frak g p \oplus S_p \cong \frak m p \oplus
\lambd p\oplus (\frak m p)^* \oplus (\lambd p)^*
 \oplus (S^*_{\qq})^v\oplus
S_{\qq}^h,
$$
and the canonical map from $\frak g = \frak m \oplus \lambd$ to
$\frak m p \oplus \lambd p$ is an isomorphism. The explicit
expression (4.8.3) for the derivative of the momentum mapping now
implies the claim at once. It also implies the claim in the
general case where $G_{\qq}$ is non-trivial since in this case the
$G$-momentum mapping $\mu_{[0]}$ arises from the $G$-momentum
mapping for the free case by $G_{\qq}$-reduction. \qed \enddemo

Let $\omega_{\qq}$ be the cotangent bundle symplectic structure on
$\sigma_{\qq} \oplus \sigma_{\qq}^*$. The decomposition (4.7.2) of
the tangent space
$$
\roman T_p(\QQ^{\roman T^*}) \cong \roman
T_{[e,\nu,0,p_{\qq}]}(G^{\roman T^*}\times S^{\roman T^*}_{\qq})_0
$$
may plainly as well be written as
$$
 \frak g p \oplus S_p \cong \frak g p\oplus (\frak g p)^*
 \oplus \sigma^*_{\qq}\oplus
\sigma_{\qq}. \tag4.9.1
$$
Given $X\in \frak g$, we will write its image in $\frak g p$ under
the canonical projection as $X_p$ and, likewise, given $\phi \in
\frak g^*$, we will denote by $\phi_p$ the unique element of
$(\frak g p)^*$ which, under the canonical injection from $(\frak
g p)^*$ into $\frak g^*$, goes to $\phi$.

\proclaim{Proposition 4.9} In terms of the right-hand side of the
decomposition {\rm (4.9.1)}, when $\QQ^{\roman T^*}$ is identified
with $\left(G^{\roman T^*}\times S^{\roman T^*}_{\qq}\right)_0$
via {\rm (4.8.2)}, at the point $ [e,\nu,0,p_{\qq}]$ of
$\left(G^{\roman T^*}\times S^{\roman T^*}_{\qq}\right)_0, $ which
corresponds to the point $p$ of $\QQ^{\roman T^*}$, the symplectic
structure
$$
\omega\colon \roman T_{[e,\nu,0,p_{\qq}]}(G^{\roman T^*}\times
S^{\roman T^*}_{\qq})_0 \times \roman
T_{[e,\nu,0,p_{\qq}]}(G^{\roman T^*}\times S^{\roman T^*}_{\qq})_0
@>>> \Bbb R
$$
is given by the formula
$$
\omega ((Y_p, \psi_p,u),(Z_p,\phi_p,v)) = \phi(Y) - \psi(Z)+
[Z,\alpha](Y) + \omega_{\qq}(u,v) \tag4.9.2
$$
where $Y,Z \in \frak g$, $\phi,\psi \in \frak g^*$, and $u,v \in
\sigma_{\qq}\oplus \sigma^*_{\qq}$,  and where, with a slight
abuse of notation, we do not distinguish in notation between $p$
and $[e,\nu,0,p_{\qq}]$.
\endproclaim

\demo{Proof} At the point $\alpha=\alpha_e\in \roman T^*_eG =\frak
g^*$, the cotangent bundle symplectic structure
$$
\omega\colon \left(\roman T_eG \times \roman T_{\alpha}\frak
g^*\right) \oplus \left(\roman T_eG \times \roman T_{\alpha}\frak
g^*\right)@>>> \Bbb R
$$
on $G^{\roman T^*}$  has the form
$$
\omega ((Y, \psi),( Z,\phi)) = \phi(Y) - \psi(Z)+
[Z,\alpha](Y),\quad Y,Z \in \frak g,\, \phi,\psi \in \frak g^*,
\tag4.9.3
$$
and is actually determined by this expression since $\omega$ is
$G$-biinvariant. This implies the assertion since $\left(G^{\roman
T^*}\times S^{\roman T^*}_{\qq}\right)_0$ arises from $G^{\roman
T^*}\times S^{\roman T^*}_{\qq}$ by $G_{\qq}$-reduction. \qed
\enddemo

\noindent {\smc 4.10. The symplectic slice.} By virtue of Lemma
4.8, restricted to the subspace $\lambd p \oplus (\lambd p)^*$ of
$$
\roman T_p(G\times_{G_p} S_p) \cong
 \frak g p \oplus S_p \cong \frak m p \oplus
\lambd p\oplus (\frak m p)^* \oplus (\lambd p)^*
 \oplus \sigma^*_{\qq}\oplus
\sigma_{\qq},
$$
the derivative of the momentum mapping comes down to the linear
map
$$
\lambd p \oplus (\lambd p)^* @>>>  (\lambd p)^* \subseteq \frak
g^*, \quad (X_p, \alpha_p) \longmapsto [X,\nu] + \alpha,
\tag4.10.1
$$
where $X\in \lambd$ and $\alpha\in \lambd^*$. The kernel $K$ of
(4.10.1) consists of the pairs $(-X_p,[X,\nu]_p)$ with $X\in
\lambd$. Proposition 4.9 implies that, given $(-X_p,[X,\nu]_p)$
and $(-Y_p,[Y,\nu]_p)$ in $K$,
$$
\omega ((-X_p,[X,\nu]_p),(-Y_p,[Y,\nu]_p)) = \nu ([Y,X]),
$$
that is, under the projection from $\lambd p \oplus (\lambd p)^*$
to the second summand, $K$ is mapped isomorphically onto $(\lambd
p)^*$ in such a way that the symplectic form on $\roman
T_p(G\times_{G_p} S_p)$, restricted to $K$, is identified with the
negative Kostant-Kirillov-Souriau form. This amounts of course to
the familiar fact that reduction at the value $\nu \in \frak g^*$
for the $G$-action on $\roman T^*G$ by left translation yields the
coadjoint orbit $G\nu$ endowed with the negative
Kostant-Kirillov-Souriau form.

It is now immediate that, in view of the decomposition (4.7.2) of
the tangent space to $\QQ^{\roman T^*}$ at $p$, the vector space
$$
V = K \oplus \sigma^*_{\qq}\oplus \sigma_{\qq} \tag4.10.2
$$
is a {\it symplectic slice\/} at the point $p$, where $K$ is
endowed with the negative Kostant-Kirillov-Souriau form or what
corresponds to it, and where $\sigma^*_{\qq}\oplus \sigma_{\qq}$
carries the linear cotangent bundle symplectic structure written
above as $\omega_{\qq}$. That is to say: Endowed with the
corresponding $G_p$-momentum mapping, $G_p$-reduction, applied to
$V$, yields a local model for the $G_{\nu}$-reduced space at
$\nu$, as a stratified symplectic space. This symplectic slice
admits a characterization in terms of the $G$-action on the base
$\QQ$ and of the geometry of the coadjoint orbit generated by
$\nu$, as explained in the introduction.

\noindent {\smc 4.11. The Witt-Artin decomposition.} Letting
$W=(\frak m p)^*$ ($=(\frak g_{\nu} p)^* $), we obtain the {\it
Witt-Artin decomposition\/}
$$
\roman T_p(\roman T^* \QQ) \cong \frak g_{\nu} p \oplus \lambd p
\oplus V \oplus W \tag4.11.1
$$
at the point $p$ of $\roman T^* \QQ$. Indeed, in view of (4.9.2),
with reference to the decomposition (4.9.1) of $\roman T_p(\roman
T^* \QQ)$, the skew-orthogonal complement $(\frak gp)^{\omega}$ of
$\frak gp$ in $\roman T_p(\roman T^* \QQ)$ consists of the triples
$$
(X_p,\psi_p, u) \in \frak g p \oplus (\frak g p)^*\oplus
(\sigma^*_{\qq}\oplus \sigma_{\qq})
$$
satisfying the equation
$$
\psi(Y) + [X,\nu](Y) =0
$$
for every $Y \in \frak g$ whence
$$
(\frak gp)^{\omega} = \{(-X_p,[X,\nu]_p, u); X \in \frak g,\, u
\in \sigma^*_{\qq}\oplus \sigma_{\qq} \} = \frak mp \oplus K\oplus
\sigma^*_{\qq}\oplus \sigma_{\qq}= \frak mp \oplus V,
$$
and $V$ is a $G_p$-invariant complement of $\frak g_{\nu}p = \frak
mp$ in $(\frak gp)^{\omega}$. Likewise, with reference to the
decomposition (4.9.1) of $\roman T_p(\roman T^* \QQ)$, in view of
(4.9.2), the skew-orthogonal complement ${(V\oplus \lambd
p)^{\omega}}$ of $V\oplus \lambd p$ in $\roman T_p(\roman T^*
\QQ)$ consists of the triples
$$
(X_p,\psi_p, 0) \in \frak g p \oplus (\frak g p)^*\oplus
(\sigma^*_{\qq}\oplus \sigma_{\qq})
$$
satisfying the equation
$$
\phi(X) -\psi(Z) - [X,\nu](Z)  =  0
$$
for every $Z \in \frak q$ and every $\phi \in \frak q^*$, and a
little thought reveals that, therefore,
 $(V\oplus
\lambd p)^{\omega} = W \oplus \frak g_{\nu}p$; moreover, $W$ is a
$G_p$-invariant Lagrangian complement of $ \frak g_{\nu}p$ in
$(V\oplus \lambd p)^{\omega}$. Hence (4.11.1) is indeed a
Witt-Artin decomposition, cf. \cite\ortrathr.

 The Witt-Artin decomposition
(4.11.1) can be characterized entirely in terms of the $G$-action
on the base $\QQ$ and the geometry of the coadjoint orbit
generated by $\nu$. Indeed, let $G_{\nu,\qq}$ ($= G_{\nu} \cap
G_\qq$) denote the stabilizer of $\qq$ for the $G_\nu$-action. The
tangent space $\frak g_\nu p$ to the $G_\nu$-orbit of $p$ at the
point $p$ decomposes as $\frak g_{\nu,\qq} p \oplus (\frak
g_{\nu,\qq} p)^{\bot}$. Here  $\frak g_{\nu,\qq} p$ is the tangent
space at $p$ to the $G_{\nu,\qq}$-orbit of $p$ in $\roman
T_\qq^*\QQ$, $\frak g_{\nu,\qq} p$ being viewed as a linear
subspace of the tangent space $\roman T_{p}\roman T_\qq^*\QQ$,
that is, of $\roman T_\qq^*\QQ$, and $(\frak g_{\nu,\qq}
p)^{\bot}$ is the orthogonal complement of $\frak g_{\nu,\qq} p$;
under the cotangent bundle projection map, $\frak g_{\nu,\qq} p$
is identified with the tangent space $\frak g_{\nu} \qq$
($\subseteq \roman T_\qq\QQ$) to the $G_{\nu}$-orbit at the point
$\qq$ of the base $\QQ$. Likewise, let $\lambd_{\nu,\qq}$ be the
orthogonal complement of $\frak g_{\nu,\qq}$ in $\frak g_\qq$, so
that $\frak g_\qq = \frak g_{\nu,\qq} \oplus \lambd_{\nu,\qq}$.
Then the vector space $\lambd p$, viewed as a linear subspace of
$\frak g p$ and hence as a linear subspace of $\roman T_\qq^*\QQ$,
decomposes as $\lambd_{\nu,\qq} p \oplus (\lambd_{\nu,\qq}
p)^{\bot}$ and, under the cotangent bundle projection map,
$(\lambd_{\nu,\qq} p)^{\bot}$ is identified with $\lambd \qq$,
viewed as a subspace of $\roman T_b\QQ$. Consequently, as
$G_{\nu,\qq}$-representations, $\frak g_\nu p \cong \frak
g_{\nu,\qq} p  \oplus \frak g_{\nu} \qq$ and ${\lambd p \cong
\lambd_{\nu,\qq} p  \oplus \lambd \qq}$.

\medskip

\centerline{\smc References} \smallskip \widestnumber\key{999}
\comment
 \ref \no  \abramars \by R. Abraham and J. E. Marsden
\book Foundations of Mechanics \publ Benjamin-\linebreak Cummings
Publishing Company
\yr 1978
\endref

 \ref \no \domboone \by P. Dombrowski \paper On the
geometry of the tangent bundle \jour J. reine angew. Mathematik
\vol 210 \yr 1962 \pages 73--88
\endref
\endcomment

\ref \no \guisteft \by V. W. Guillemin and S. Sternberg \paper A
normal form for the Moment Map \paperinfo in: Differential
Geometric Methods in Mathematical Physics, S. Sternberg, ed. \publ
Reidel Publishing Company \publaddr Dordrecht \yr 1984
\endref

\ref \no \marleone \by C. M. Marle \paper Mod\`ele d'action
hamiltonienne d'un groupe de Lie sur une vari\'et\'e symplectique
\jour Rendiconti del Seminario Matematico \vol 43 \yr 1985 \pages
227--251 \publ Universit\`a e Politechnico \publaddr Torino
\endref

\ref \no \ortrathr \by J.-P. Ortega and T. Ratiu \book Hamiltonian
Reduction \bookinfo Progress in Mathematics, Vol. 222 \publ
Birkh\"auser Verlag \publaddr Boston $\cdot$ Basel $\cdot$ Berlin
\yr 2004
\endref

\ref \no \perrosou
\by M. Perlmutter, M. Rodr\'\i guez-Olmos, and
M. E. Sousa-Dias \paper The Witt-Artin decomposition of a
cotangent-lifted action \paperinfo {\tt math.SG/0501207}
\endref

\ref \no \schmah
\by T. Schmah \paper A cotangent bundle slice
theorem \paperinfo {\tt math.SG/0409148}
\endref

\ref \no \sjamlerm \by R. Sjamaar and E. Lerman \paper Stratified
symplectic spaces and reduction \jour Ann. of Math. \vol 134 \yr
1991 \pages 375--422
\endref

\enddocument